\documentclass[12pt,reqno]{amsart}%
\usepackage{amsfonts,anysize}
\usepackage[charter]{mathdesign}
\usepackage[colorlinks,urlcolor=blue,citecolor=blue,linkcolor=blue]{hyperref}

\marginsize{2.5cm}{2.5cm}{1cm}{1cm}
\newtheorem{thm}{Theorem}[section]

\newtheorem{lem}[thm]{Lemma}

\newtheorem{prop}[thm]{Proposition}
\theoremstyle{definition}

\theoremstyle{remark}
\newtheorem{rem}[thm]{Remark}
\newtheorem*{pf}{Proof}
\numberwithin{equation}{section}
\allowdisplaybreaks

\begin{document}

\title[Standing waves for quasilinear Schr\"{o}inger equations with indefinite potentials]{\large Standing waves for quasilinear Schr\"{o}inger equations with indefinite potentials}
\thanks{This work was supported by NSFC (11671331) and NSFFJ (2014J06002).}
\author{Shibo Liu
\and Jian Zhou}
\dedicatory{School of Mathematical Sciences, Xiamen University, Xiamen 361005, China}

\begin{abstract}
We consider quasilinear Schr\"{o}dinger equations in $\mathbb{R}^{N}$ of the
form%
\[
-\Delta u+V(x)u-u\Delta(u^{2})=g(u)\text{,}%
\]
where $g(u)$ is $4$-superlinear. Unlike all known results in the literature,
the Schr\"{o}dinger operator $-\Delta+V$ is allowed to be indefinite, hence the
variational functional does not satisfy the mountain pass geometry. By a local
linking argument and Morse theory, we obtain a nontrivial solution for the
problem. In case that $g$ is odd, we get an unbounded sequence of solutions.

\end{abstract}
\maketitle

\section{Introduction}

In this paper we consider quasilinear stationary Schr\"{o}dinger equations in
$\mathbb{R}^{N}$ of the form%
\begin{equation}
-\Delta u+V(x)u-u\Delta(\left\vert u\right\vert ^{2})=g(u)\text{.} \label{e1}%
\end{equation}
This kind of equations arise when we are looking for standing waves
$\psi(t,x)=\mathrm{e}^{-i\omega t}u(x)$ for the time-dependent quasilinear
Schr\"{o}dinger equation%
\begin{equation}
i\partial_{t}\psi=-\Delta\psi+U(x)\psi-\psi\Delta(\vert \psi\vert
^{2})-\bar{g}(\vert \psi\vert ^{2})\psi\text{,\qquad}\left(
t,x\right)  \in\mathbb{R}^{+}\times\mathbb{R}^{N}\text{,} \label{sd}%
\end{equation}
which was used for the superfluid film equation in plasma physics by Kurihar
\cite{ref4}. As mentioned by Ruiz-Siciliano \cite{MR2630099}, it is also
introduced in \cite{ref7,MR1986307,ref2} to study a model of self-trapped
electrons in quadratic or hexagonal lattices. \ Note that the potential $V$ in
\eqref{e1} is given by
\begin{equation}
V(x)=U(x)-\omega\text{.} \label{vv}%
\end{equation}
Thus, if the frequency $\omega$ is large, $V$ may be indefinite in sign. For
physical reason, it is natural to make the following assumption on the
nonlinearity
\begin{equation}
\lim_{t\rightarrow0}\frac{g(t)}{t}=0\text{.} \label{e2}%
\end{equation}

To the best of our knowledge, for nonlinearity of the form $g(u)=\left\vert
u\right\vert ^{p-2}u$, the first mathematical studies of the equation \eqref{e1}
seems to be Poppenberg \emph{et.\ al.}\ \cite{MR1899450} for the one
dimensional case and Liu-Wang \cite{MR1933335} for higher dimensional case.
The proofs in these papers are based on constrained minimization argument.

To consider more general nonlinearities, we should seek for a free variational
formulation for the problem. Formally, the solutions of the problem \eqref{e1}
should be critical points of the following functional
\[
J(u)=\frac{1}{2}\int\left(  1+2u^{2}\right)  \left\vert \nabla u\right\vert
^{2}+\frac{1}{2}\int V(x)u^{2}-\int G(u)\text{,\qquad where }G(t)=\int_{0}%
^{t}g(\tau)d\tau\text{.}%
\]
Unfortunately, the functional $J$ could not be defined for every $u\in
H^{1}(\mathbb{R}^{N})$. Therefore the standard variational methods could not
be applied. To overcome this difficulty, Liu \emph{et.\ al.}\ \cite{MR1949452}
and Colin-Jeanjean \cite{MR2029068} introduced a nonlinear transformation $f$
so that if $v\in H^{1}(\mathbb{R}^{N})$ is a critical point of $\Phi
:H^{1}(\mathbb{R}^{N})\rightarrow\mathbb{R}$,%
\begin{equation}
\Phi(v)=\frac{1}{2}\int\left\vert \nabla v\right\vert ^{2}+\frac{1}{2}\int
V(x)f^{2}(v)-\int G(f(v))\text{,} \label{qq}%
\end{equation}
then $u=f(v)$ is a solution of \eqref{e1}.

Since the publication of \cite{MR2029068,MR1949452}, many results about
\eqref{e1} have been obtained by various authors using this transformation. For
example, the case that the potential $V$ is $\mathbb{Z}^{N}$-periodic (or
asymptotically periodic) and $g$ is $4$-superlinear is studied in Silva-Vieira
\cite{MR2580150}, Yang \cite{MR2927594}, Fang-Szulkin \cite{MR3003301} and
Zhang \emph{et.\ al.}\ \cite{MR3620755}. In \cite{MR3160456}, Wu studied the
case that the nonlinearity is $4$-superlinear and the potential has a positive
lower bound and satisfies the condition $\left(  V\right)  $ below. In Furtado
\emph{et.\ al.}\ \cite{MR3327366}, the asymptotically linear case (the
nonlinearity behaves like $t$ at the origin and like $t^{3}$ at infinity) is
investigated, where the potential $V$ satisfies the same conditions as in
\cite{MR3160456}. For equations with concave and convex nonlinearities, one
can consult do \'{O} and Severo \cite{MR2461565}. For problems with critical
nonlinearities, see Silva-Vieira \cite{MR2659677} and Wang-Zou \cite{MR2885550}%
. For the supercritical case, we refer the reader to Moameni \cite{MR2358356}
and Miyagaki-Moreira \cite{MR3250499}.

In all these papers, it is required that the potential $V$ satisfies the
positive condition%
\begin{equation}
\alpha:=\inf_{\mathbb{R}^{N}}V>0\text{.} \label{v0}%
\end{equation}
With this condition and suitable conditions on the nonlinearity, one can show
that the zero function $v=0$ is a local minimizer of the functional $\Phi$,
and $\Phi$ would then verify the mountain pass geometry and the mountain pass
theorem can be applied to produce a solution. However, from \eqref{vv} we can
see that, if we want to find standing waves $\psi(t,x)=\mathrm{e}^{-i\omega
t}u(x)$ of \eqref{sd} with large $\omega$, then the potential $V$ in \eqref{e1}
could not satisfy \eqref{v0}.

In the literature there are some existence results which allow the potential
$V$ to be negative somewhere. The strategy is to write $V=V^{+}-V^{-}$ with
$V^{\pm}=\max\left\{  0,\pm V\right\}  $. Then if $V^{-}$ is in some sense
small, it can be absorbed and the functional $\Phi$ still verifies the
mountain pass geometry. We refer the reader to Fang-Han \cite{MR3159414} and
Maia \emph{et.\ al.} \cite{MR3473792} for this kind of results.

In Zhang \emph{et.\ al.} \cite{MR3240105}, the authors studied the problem
\eqref{e1} with sign-changing potential. Their potential satisfies conditions
slightly general than our condition $\left(  V\right)  $ below. Because $V$ is
sign-changing, the function%
\[
u\mapsto\left(  \int\left(  \left\vert \nabla u\right\vert ^{2}+V(x)u^{2}%
\right)  \right)  ^{1/2}%
\]
is no longer a norm on the function space. This will bring some problems for
verifying the boundedness of asymptotically critical sequences. To overcome
this difficulty, they chose a constant $V_{0}>0$ such that
\[
\tilde{V}(x)=V(x)+V_{0}>0
\]
and consider the equivalent problem%
\[
-\Delta u+\tilde{V}(x)u-u\Delta(\left\vert u\right\vert ^{2})=\tilde
{g}(x,u)\text{,\qquad}x\in\mathbb{R}^{N}\text{,}%
\]
where $\tilde{g}(x,u)=g(x,u)+V_{0}u$. Unfortunately, from%
\[
\mathcal{G}(x,u):=\frac{1}{4}g(x,u)u-G(x,u)\geq0
\]
in their condition $\left(  G_{2}\right)  $, we could not deduce%
\[
\mathcal{\tilde{G}}(x,u):=\frac{1}{4}\tilde{g}(x,u)u-\tilde{G}%
(x,u)=\mathcal{G}(x,u)-\frac{1}{4}V_{0}u^{2}\geq0\text{.}%
\]
Therefore unlike what the authors declared at the beginning of \cite[\S 2]%
{MR3240105}, this new nonlinearity $\tilde{g}(x,u)$ does not satisfy their
condition $\left(  G_{2}\right)  $ any more. Because $\mathcal{\tilde{G}%
}(x,u)\geq0$ is crucial to get%
\[
\mathrm{meas}(\Omega_{n}(r,+\infty))\rightarrow0
\]
in \cite[Page 1768]{MR3240105}, their result can only be valid for the case
that $V$ is positive.

In conclusion, for the quasilinear Schr\"{o}dinger equation \eqref{e1}, as far
as we know, up to now in the literature there is no research devoted to the
situation that the zero function $v=0$ fails to be a local minimizer of $\Phi$.

The purpose of this paper is to present the first results in this indefinite
situation. Firstly, we present our assumptions on the potential $V(x)$ and the
nonlinearity $g(u)$.

\begin{itemize}
\item[$(V)$] $V\in C(\mathbb{R}^{N})$ is bounded from below and, $\mathfrak{m}
(V^{-1}(-\infty,M])<\infty$ for all $M>0$, where $\mathfrak{m}$ is the
Lebesgue measure on $\mathbb{R}^{N}$.

\item[$\left(  g_{0}\right)  $] $g\in C(\mathbb{R})$ and there exist $C>0$ and
$p\in\left(  4,2\cdot2^{\ast}\right)  $ such that%
\begin{equation}
\vert g(t)\vert \leq C\left(  \vert t\vert +\vert
t\vert ^{p-1}\right)  \text{,} \label{e3}%
\end{equation}
where the critical Sobolev exponent $2^{*}=2N/(N-2)$ for $N\ge3$ and
$2^{*}=\infty$ for $N=2$.

\item[$\left(  g_{1}\right)  $] there exists $\mu>4$ such that for $t\neq0$
there holds%
\[
0<\mu G(t)\leq tg(t)\text{.}%
\]

\end{itemize}

As we shall see, for problem \eqref{e1} with indefinite potential $V$, there are
several interesting features. Firstly, it is well known that the mountain pass
theorem is not applicable in this situation. For semilinear problems, we may
try to apply the linking theorem. But for our quasilinear problem, because%
\[
Q(v)=\frac{1}{2}\int\left\vert \nabla v\right\vert ^{2}+\frac{1}{2}\int
V(x)f^{2}(v)\text{,}%
\]
the principle part of $\Phi$, is not a quadratic form on $v$, the linking
theorem is also not applicable. The reason is that for applying the linking
theorem, one needs to decompose the function space according to certain
quadratic form, while in our functional $\Phi$, there is no natural quadratic
form. Even if we decompose the space according to some quadratic forms like
$\mathfrak{B}$ given by%
\[
\mathfrak{B}(u)=\frac{1}{2}\int\left(  \left\vert \nabla u\right\vert
^{2}+V(x)u^{2}\right)  \text{,}%
\]
the non-quadratic feature of $Q$ will also prevent us to verify the linking
geometry. Our key observation is that, just like our previous study
\cite{MR3303004} on the Schr\"{o}dinger-Poisson systems%
\begin{equation}
\left\{
\begin{array}
[c]{ll}%
-\Delta u+V(x)u+\phi u=g(u) & \text{in }\mathbb{R}^{3}\text{,}\\
-\Delta\phi=u^{2} & \text{in }\mathbb{R}^{3}\text{,}%
\end{array}
\right. \label{sp}%
\end{equation}
the functional $\Phi$ has a local linking \cite{MR1312028,MR802575} at the
origin with respect to the decomposition $X=X^{-}\oplus X^{+}$, where $X$ is
our working space that will introduced later, $X^{-}$ and $X^{+}$ are the
negative and positive spaces of the quadratic form $\mathfrak{B}$, respectively.

Secondly, for indefinite semilinear Schr\"{o}dinger equations
\[
-\Delta u+V(x)u=g(u)\qquad\text{in $\mathbb{R}^{N}$}
\]
and the Schr\"{o}dinger-Poisson systems \eqref{sp}, the principle part of the
variational functional is the quadratic form $\mathfrak{B}$. Therefore, to
obtain the boundedness of asymptotically critical sequences $\left\{
u_{n}\right\}  $, one usually needs to decompose $u_{n}=u_{n}^{-}+u_{n}^{+}$,
where $u_{n}^{\pm}$ is the orthogonal projection of $u_{n}$ on $X^{\pm}$, see
e.g., \cite{MR1751952,MR2957647,MR3656292}. In our quasilinear case, again,
the nonquadratic feature of $Q$ makes such decomposition useless for verifying
the boundedness of the sequences.

Before state our results, we shall introduce the suitable function space in
which we will find critical points of $\Phi$. Since $V$ is bounded from below,
as in \cite[Page 45]{MR3303004} we choose $m>0$ such that%
\begin{equation}
\tilde{V}(x)=V(x)+m>1 \label{V}%
\end{equation}
for all $x\in\mathbb{R}^{N}$. On the subspace%
\[
X=\left\{  u\in H^{1}(\mathbb{R}^{N})\left\vert \,\int V(x)u^{2}%
<\infty\right.  \right\}
\]
we equip the inner product%
\[
\left\langle u,v\right\rangle =\int\left(  \nabla u\cdot\nabla v+\tilde
{V}(x)uv\right)
\]
and corresponding norm $\Vert u\Vert =\left\langle u,u\right\rangle
^{1/2}$. Then $X$ is a Hilbert space and by Bartsch-Wang \cite{MR1349229} we
have a compact embedding $X\hookrightarrow L^{s}(\mathbb{R}^{N})$ for
$s\in\lbrack2,2^{\ast})$.

By the compactness of the embedding $X\hookrightarrow L^{2}(\mathbb{R}^{N})$,
applying the spectral theory of self-adjoint compact operators, we see that
the eigenvalue problem%
\begin{equation}
-\Delta u+V(x)u=\lambda u\text{,\qquad}u\in X \label{eg}%
\end{equation}
possesses a complete sequence of eigenvalues
\[
-\infty<\lambda_{1}\leq\lambda_{2}\leq\cdots\text{,\qquad}\lambda
_{i}\rightarrow+\infty\text{.}%
\]
Each $\lambda_{i}$ has been repeated in the sequence according to its finite
multiplicity. We denote by $\phi_{i}$ the eigenfunction of $\lambda_{i}$ with
$\vert \phi_{i}\vert _{2}=1$, where $\left\vert \,\cdot
\,\right\vert _{q}$ is the $L^{q}$-norm. Now we can state our main results.

\begin{thm}
\label{t1}Suppose that $\left(  V\right)  $, $\left(  g_{0}\right)  $,
$\left(  g_{1}\right)  $ and \eqref{e2} hold. If $0$ is not an eigenvalue of
\eqref{eg}, then \eqref{e1} has a nontrivial solution.
\end{thm}

\begin{thm}
\label{t2}Suppose that $\left(  V\right)  $, $\left(  g_{0}\right)  $ and
$\left(  g_{1}\right)  $ hold. If $g$ is odd, then \eqref{e1} has a sequence of
solutions $\left\{  u_{n}\right\}  $ such that $J(u_{n})\rightarrow+\infty$.
\end{thm}

\begin{rem}
Note that in Theorem \ref{t2} we don't need the local condition \eqref{e2}.
\end{rem}

The paper is organized as follow. In Section 2 we recall the transformation
$f$ introduced in \cite{MR2029068}, which converts the problem of solving
\eqref{e1} into searching for critical points of the functional $\Phi$ given in
\eqref{qq}. Then, as the first step of dealing with $\Phi$ we show that it
satisfies the Cerami condition. We will prove Theorem \ref{t1} by applying
Morse theory. Therefore in Section 3 we recall some necessary concepts and
results in Morse theory. Then we verify the local linking property for $\Phi$
and compute the critical groups of $\Phi$ at infinity. With these preparation
we give the proof of Theorem \ref{t1}. Finally, in Section 4, we apply the
symmetric mountain pass theorem of Ambrosetti-Rabinowitz \cite{MR0370183} to
prove Theorem \ref{t2}.

\section{Variational reformulation and Cerami condition}

Following Colin-Jeanjean \cite{MR2029068}, we make the change of variables by
$u=f(v)$, where $f$ is an odd function defined by
\[
\dot{f}(t)=\frac{1}{\sqrt{1+2f^{2}(t)}}\text{,\qquad}f(0)=0
\]
on $[0,+\infty)$. Then we have the following proposition, whose proof can be
found in \cite{MR3121525,MR2029068}.

\begin{prop}
\label{p1}The function $f$ satisfies the following properties:

\begin{enumerate}
\item $f\in C^{\infty}(\mathbb{R})$ is strictly increasing, therefore is invertible.

\item \label{ite0}$\vert f(t)\vert \leq t$ and $|\dot{f}(t)|\leq1$
for all $t\in\mathbb{R}$. Moreover,%
\[
\dot{f}(0)=\lim_{t\rightarrow0}\frac{f(t)}{t}=1\text{.}%
\]

\item \label{ite1}for all $t>0$ we have%
\[
\frac{1}{2}f(t)\leq\dot{f}(t)t\leq f(t)\text{.}%
\]

\item \label{ite2}for all $t\in\mathbb{R}$ we have $f^{2}(t)\geq f(t)\dot
{f}(t)t$ and $\vert f(t)\vert \leq2^{1/4}\left\vert t\right\vert
^{1/2}$.

\item \label{ite3}there exists a positive constant $\kappa$ such that%
\[
\vert f(t)\vert \geq\kappa\left\vert t\right\vert \quad\text{for
}\left\vert t\right\vert \leq1\text{,\qquad\qquad}\vert f(t)\vert
\geq\kappa\left\vert t\right\vert ^{1/2}\text{ for }\left\vert t\right\vert
\geq1\text{.}%
\]

\item \label{ite4}for each $\lambda>0$, there exists a positive constant
$C_{\lambda}$ such that $f^{2}(\lambda t)\leq C_{\lambda}f^{2}(t)$.
\end{enumerate}
\end{prop}

By the growth condition of the nonlinearity $g$ and the properties of the
transformation $f$, it is easy to verify that the functional
\[
\Phi(v)=J(f(v))=\frac{1}{2}\int\left\vert \nabla v\right\vert ^{2}+\frac{1}%
{2}\int V(x)f^{2}(v)-\int G(f(v))
\]
is well defined and of class $C^{1}$ on the Sobolev space $X$ introduced in
Section 1, and if $v$ is a critical point of $\Phi$, then $u=f(v)$ is a
solution of our problem \eqref{e1}, see \cite{MR2029068} for the details.

Hence, to prove our main results, we shall look for critical points of the
functional $\Phi$. Firstly, we need to show that the functional $\Phi$
satisfies the Cerami condition. To this end, we set%
\[
\tilde{g}(t)=g(t)+mt\text{,\qquad}\tilde{G}(t)=\int_{0}^{t}\tilde{g}%
(\tau)d\tau=G(t)+\frac{m}{2}t^{2}%
\]
and rewite $\Phi$ in the following form%
\[
\Phi(v)=\frac{1}{2}\int\left\vert \nabla v\right\vert ^{2}+\frac{1}{2}%
\int\tilde{V}(x)f^{2}(v)-\int\tilde{G}(f(v))\text{,}%
\]
where $\tilde{V}$ is given in \eqref{V}. Note that by $\left(  g_{1}\right)  $,
the new nonlinearity $\tilde{g}$ satisfies%
\begin{equation}
\tilde{G}(t)-\frac{1}{\mu}\tilde{g}(t)t\leq\left(  \frac{1}{2}-\frac{1}{\mu
}\right)  mt^{2}\text{.} \label{ge}%
\end{equation}

\begin{lem}
\label{l1}Under the assumptions $\left(  V\right)  $, $\left(  g_{0}\right)  $
and $\left(  g_{1}\right)  $, the functional $\Phi$ satisfies the Cerami condition.
\end{lem}

\begin{pf}
Let $\left\{  v_{n}\right\}  $ be a Cerami sequence of $\Phi$, that is,
\[
\Phi(v_{n})\rightarrow c\text{,\qquad}\left(  1+\Vert v_{n}\Vert
\right)  \Phi^{\prime}(v_{n})\rightarrow0
\]
for some $c\in\mathbb{R}$. We claim that there exists $C>0$ such that%
\begin{equation}
\rho_{n}:=\left\{  \int\left(  \vert \nabla v_{n}\vert ^{2}%
+\tilde{V}(x)f^{2}(v_{n})\right)  \right\}  ^{1/2}\leq C\text{.} \label{e4}%
\end{equation}
If this is not true, we may assume $\rho_{n}\rightarrow+\infty$. Consider the
sequence%
\[
h_{n}=\frac{f(v_{n})}{\rho_{n}}\text{.}%
\]
Then since $|\dot{f}(t)|\leq1$, we have%
\begin{align}
\Vert h_{n}\Vert ^{2}  &  =\frac{1}{\rho_{n}^{2}}\int\left(
\vert \nabla f(v_{n})\vert ^{2}+\tilde{V}(x)f^{2}(v_{n})\right)
\nonumber\\
&  =\frac{1}{\rho_{n}^{2}}\int\left(  |\dot{f}(v_{n})|^{2}\vert \nabla
v_{n}\vert ^{2}+\tilde{V}(x)f^{2}(v_{n})\right) \nonumber\\
&  \leq\frac{1}{\rho_{n}^{2}}\int\left(  \vert \nabla v_{n}\vert
^{2}+\tilde{V}(x)f^{2}(v_{n})\right)  =1\text{.} \label{e7}%
\end{align}
Consequently, $\left\{  h_{n}\right\}  $ is bounded in $X$. Up to a
subsequence, by the compactness of the embedding $X\hookrightarrow
L^{2}(\mathbb{R}^{N})$, we may assume that%
\[
h_{n}\rightharpoonup h\text{ in }X\text{,\qquad}h_{n}\rightarrow h\text{ in
}L^{2}\text{,\qquad}h_{n}\rightarrow h\text{ a.e. in }\mathbb{R}^{N}\text{.}%
\]
As in Colin-Jeanjean \cite[Page 225]{MR2029068}, set%
\[
\phi_{n}=\sqrt{1+2f^{2}(v_{n})}f(v_{n})\text{.}%
\]
Then using Proposition \ref{p1} \eqref{ite0} and \eqref{ite2} we have%
\begin{align*}
\vert \nabla\phi_{n}\vert  &  =\left(  1+\frac{2f^{2}(v_{n}%
)}{1+2f^{2}(v_{n})}\right)  \vert \nabla v_{n}\vert \leq2\vert
\nabla v_{n}\vert \text{,}\\
\Vert \phi_{n}\Vert ^{2}  &  =\int\left(  \vert \nabla\phi
_{n}\vert ^{2}+\tilde{V}(x)\phi_{n}^{2}\right) \\
&  \leq4\int\vert \nabla v_{n}\vert ^{2}+\int\tilde{V}%
(x)f^{2}(v_{n})+2\int\tilde{V}(x)f^{4}(v_{n})\\
&  \leq4\int\vert \nabla v_{n}\vert ^{2}+\int\tilde{V}(x)v_{n}%
^{2}+2\int\tilde{V}(x)(2v_{n}^{2})\leq5\Vert v_{n}\Vert
^{2}\text{.}%
\end{align*}
Therefore, by \eqref{ge}, we have%
\begin{align*}
c+o(1)  &  =\Phi(v_{n})-\frac{1}{\mu}\langle\Phi^{\prime}(v_{n}),\phi
_{n}\rangle\\
&  =\frac{1}{2}\int\left(  \vert \nabla v_{n}\vert ^{2}+\tilde
{V}(x)f^{2}(v_{n})\right)  -\int\tilde{G}(f(v_{n}))\\
&  \qquad-\frac{1}{\mu}\left\{  \int\left(  \left(  1+\frac{2f^{2}(v_{n}%
)}{1+2f^{2}(v_{n})}\right)  \vert \nabla v_{n}\vert ^{2}+\tilde
{V}(x)f^{2}(v_{n})\right)  -\int\tilde{g}(f(v_{n}))f(v_{n})\right\} \\
&  =\int\left(  \frac{1}{2}-\frac{1}{\mu}\left(  1+\frac{2f^{2}(v_{n}%
)}{1+2f^{2}(v_{n})}\right)  \right)  \vert \nabla v_{n}\vert
^{2}+\left(  \frac{1}{2}-\frac{1}{\mu}\right)  \int\tilde{V}(x)f^{2}(v_{n})\\
&  \qquad\qquad+\int\left(  \frac{1}{\mu}\tilde{g}(f(v_{n})f(v_{n})-\tilde
{G}(f(v_{n}))\right) \\
&  \geq\left(  \frac{1}{2}-\frac{2}{\mu}\right)  \int\left(  \vert \nabla
v_{n}\vert ^{2}+\tilde{V}(x)f^{2}(v_{n})\right)  -\left(  \frac{1}%
{2}-\frac{1}{\mu}\right)  \int mf^{2}(v_{n})\text{.}%
\end{align*}
Multiplying both sides by $\rho_{n}^{-2}$, we deduce%
\[
\left(  \frac{1}{2}-\frac{1}{\mu}\right)  m\vert h_{n}\vert
_{2}^{2}=\frac{m}{\rho_{n}^{2}}\left(  \frac{1}{2}-\frac{1}{\mu}\right)  \int
f^{2}(v_{n})\geq\frac{1}{2}-\frac{2}{\mu}\text{.}%
\]
Since $h_{n}\rightarrow h$ in $L^{2}$ and $\mu>4$, it follows that
$h\neq\mathbf{0}$, the zero function. Hence the set%
\[
\Theta=\left\{  \left.  x\in\mathbb{R}^{N}\right\vert h(x)\neq0\right\}
\]
is of positive Lebesgue measure.

By our assumption $\left(  g_{1}\right)  $, it is well known that%
\[
\frac{\tilde{G}(t)}{t^{2}}=\frac{1}{t^{2}}\left(  G(t)+\frac{1}{2}%
mt^{2}\right)  \rightarrow+\infty
\]
as $\vert t\vert \rightarrow\infty$. For $x\in\Theta$, we have
$h_{n}(x)\rightarrow h(x)\neq0$ and%
\[
\vert f(v_{n}(x))\vert =\rho_{n}\vert h_{n}(x)\vert
\rightarrow\infty\text{,}\qquad\frac{\tilde{G}(f(v_{n}(x)))}{f^{2}(v_{n}%
(x))}h_{n}^{2}(x)\rightarrow+\infty\text{.}%
\]
By the Fatou lemma and noting that $\tilde{G}(t)\geq0$, we deduce%
\begin{equation}
\int\frac{\tilde{G}(f(v_{n}))}{\rho_{n}^{2}}\geq\int_{\Theta}\frac{\tilde
{G}(f(v_{n}))}{f^{2}(v_{n})}h_{n}^{2}\rightarrow+\infty\text{.} \label{e8}%
\end{equation}
Therefore for large $n$ we have%
\begin{align*}
c-1\leq\Phi(v_{n})  &  =\frac{1}{2}\int\left( \vert \nabla
v_{n}\vert ^{2}+\tilde{V}(x)f^{2}(v_{n})\right)  -\int\tilde{G}%
(f(v_{n}))\\
&  =\rho_{n}^{2}\left(  \frac{1}{2}-\int\frac{\tilde{G}(f(v_{n}))}{\rho
_{n}^{2}}\right)  \rightarrow-\infty\text{.}%
\end{align*}
This is impossible. Therefore, our claim \eqref{e4} is true.

Next, we can follow the same argument as Wu \cite[Page 2626--2628]{MR3160456}
(for the special case $\alpha=1$, see Remark \ref{rk1} below) to show that
$\left\{  v_{n}\right\}  $ is bounded in $X$, and has a convergent subsequence.

More precisely, using Proposition \ref{p1} \eqref{ite3} and \eqref{ite4} we
can show that there exists a constant $C>0$ which may be depend on the sequence $\{v_n\}$, such that%
\[
\rho_{n}^{2}\geq C\Vert v_{n}\Vert ^{2}\text{,}%
\]
this and the boundedness of $\left\{  \rho_{n}\right\}  $ imply that $\left\{
v_{n}\right\}  $ is bounded. By the compactness of the embedding
$X\hookrightarrow L^{2}(\mathbb{R}^{N})$, up to a subsequence we may assume%
\[
v_{n}\rightharpoonup v\text{ in }X\text{,\qquad}v_{n}\rightarrow v\text{ in
}L^{s}(\mathbb{R}^{N})
\]
for $s\in\lbrack2,2^{\ast})$. By the growth condition \eqref{e3} and the
properties of $f$ described in Proposition \ref{p1}, we deduce%
\[
\int\left[  \tilde{g}(f(v_{n}))\dot{f}(v_{n})-\tilde{g}(f(v))\dot
{f}(v)\right]  \left(  v_{n}-v\right)  \rightarrow0\text{.}%
\]
As in the proof of \cite[Lemma 3.11]{MR3003301} we can prove that there is a
constant $C>0$ such that%
\[
\int\vert \nabla(v_{n}-v)\vert ^{2}+\int\tilde{V}(x)\left[
f(v_{n})\dot{f}(v_{n})-f(v)\dot{f}(v)\right]  \left(  v_{n}-v\right)  \geq
C\Vert v_{n}-v\Vert^2 \text{.}%
\]
Consequently,%
\begin{align*}
C\Vert v_{n}-v\Vert ^{2}  &  \leq\int\vert \nabla
(v_{n}-v)\vert ^{2}+\int\tilde{V}(x)\left[  f(v_{n})\dot{f}%
(v_{n})-f(v)\dot{f}(v)\right]  \left(  v_{n}-v\right) \\
&  \qquad\qquad-\int\left[  \tilde{g}(f(v_{n}))\dot{f}(v_{n})-\tilde
{g}(f(v))\dot{f}(v)\right]  \left(  v_{n}-v\right)  +o(1)\\
&  =\langle\Phi^{\prime}(v_{n})-\Phi^{\prime}(v),v_{n}-v\rangle+o(1)\text{.}%
\end{align*}
Since $\left\{  v_{n}\right\}  $ is a bounded Cerami sequence and
$v_{n}\rightharpoonup v$, the right hand side goes to zero as $n\to\infty$,
and we deduce $v_{n}\rightarrow v$ in $X$.
\end{pf}

\begin{rem}
\label{rk1}In \cite{MR3160456}, Wu considered slightly general equations of
the form
\[
-\Delta u+V(x)u-\vert u\vert ^{2\alpha-2}u\Delta(\vert
u\vert ^{2\alpha})=g(x,u)\text{.}%
\]
Obviously, if $\alpha=1$ this equation reduces to our \eqref{e1}. Of course, he
only studied the positive definite case $\inf_{\mathbb{R}^{N}}V>0$.
\end{rem}

\section{Proof of Theorem \ref{t1}}

To prove Theorem \ref{t1}, we will apply the infinite dimensional Morse
theory, see, e.g., Chang \cite{MR1196690} and Mawhin-Willem \cite[Chapter
8]{MR982267}. We start by recalling some concepts and results.

Let $X$ be a Banach space, $\varphi:X\rightarrow\mathbb{R}$ be a $C^{1}%
$-functional, $u$ is an isolated critical point of $\varphi$ and
$\varphi(u)=c$. Then
\[
C_{i}(\varphi,u):=H_{i}(\varphi_{c},\varphi_{c}\backslash\{0\})\text{,\qquad
}i\in\mathbb{N}=\{0,1,2,\ldots\}\text{,}%
\]
is called the $i$-th critical group of $\varphi$ at $u$, where $\varphi
_{c}:=\varphi^{-1}(-\infty,c]$ and $H_{\ast}$ stands for the singular homology
with coefficients in $\mathbb{Z}$.

If $\varphi$ satisfies Cerami condition and the critical values of $\varphi$
are bounded from below by $\alpha$, then following Bartsch-Li \cite{MR1420790}%
, we define the $i$-th critical group of $\varphi$ at infinity by
\[
C_{i}(\varphi,\infty):=H_{i}(X,\varphi_{\alpha})\text{,\qquad}i\in
\mathbb{N}\text{.}%
\]
It is well known that the homology on the right side does not depend on the choice
of $\alpha$.

\begin{prop}
[{\cite[Proposition 3.6]{MR1420790}}]\label{ap1}If $\varphi\in C^{1}%
(X,\mathbb{R})$ satisfies the Cerami condition and $C_{\ell}(\varphi,0)\neq
C_{\ell}(\varphi,\infty)$ for some $\ell\in\mathbb{N}$, then $\varphi$ has a
nonzero critical point.
\end{prop}

\begin{prop}
[{\cite[Theorem 2.1]{MR1110119}}]\label{ap2}Suppose $\varphi\in C^{1}%
(X,\mathbb{R})$ has a \emph{local linking} at $0$ with respect to the
decomposition $X=\allowbreak Y\oplus Z$, i.e., for some $\varepsilon>0$,
\[%
\begin{array}
[c]{ll}%
\varphi(u)\leq0 & \text{for }u\in Y\cap B_{\varepsilon}\text{,}\\
\varphi(u)>0 & \text{for }u\in(Z\backslash\{0\})\cap B_{\varepsilon}\text{,}%
\end{array}
\]
where $B_{\varepsilon}=\left\{  \left.  u\in X\right\vert \,\Vert u\Vert
\leq\varepsilon\right\}  $. If $\ell=\dim Y<\infty$, then $C_{\ell}%
(\varphi,0)\neq0$.
\end{prop}

Now we are ready to prove Theorem \ref{t1}. Set $\lambda_{0}=-\infty$, since
$0$ is not an eigenvalue of \eqref{eg}, we may assume that for some $\ell\geq0$
we have $0\in\left(  \lambda_{\ell},\lambda_{\ell+1}\right)  $. For $\ell\ge1$
we set%
\[
X^{-}=\operatorname*{span}\left\{  \phi_{1},\ldots,\phi_{\ell}\right\}
\text{,\qquad}X^{+}=\left(  X^{-}\right)  ^{\bot}\text{.}%
\]
If $\ell=0$, we set $X^{-}=\{0\}$ and $X^{+}=X$. Then $X^{-}$ and $X^{+}$ are
the negative space and positive space of the quadratic form%
\[
\mathfrak{B}(v)=\frac{1}{2}\int\left( \vert \nabla v\vert
^{2}+V(x)v^{2}\right)
\]
respectively, note that $\dim X^{-}=\ell$. Moreover, there is a constant
$\eta>0$ such that%
\begin{equation}
\pm\mathfrak{B}(v)\geq\eta\Vert v\Vert ^{2}\text{,\qquad}v\in
X^{\pm}\text{.} \label{qe}%
\end{equation}

\begin{lem}
\label{l2}The functional $\Phi$ has a local linking at $0$ with respect to
decomposition $X=X^{-}\oplus X^{+}$.
\end{lem}

\begin{pf}
Unlike in semilinear problems, because the principle part of $\Phi$, denoted
by $Q:X\rightarrow\mathbb{R}$,%
\[
Q(v)=\frac{1}{2}\int\left( \vert \nabla v\vert ^{2}+V(x)f^{2}%
(v)\right)  \text{,}%
\]
is not a quadratic form on $v$, it seems difficult to verify the local linking
property of $\Phi$. Fortunately, by the nice properties of the map $f$ given
by Proposition \ref{p1}, we can easily overcome this difficulty.

By direct computation, we find that
\[
\ddot{f}(t)=-\frac{f(t)\dot{f}(t)}{\left(  1+2f^{2}(t)\right)  ^{3/2}}\text{.}%
\]
By Proposition \ref{p1}, it is easy to see that $\ddot{f}$ is bounded. Hence,
although $\Phi$ is only of class $C^{1}$, its principle part $Q$ is a $C^{2}%
$-functional on $X$, with derivatives given by%
\begin{align*}
\langle Q^{\prime}(v),\phi\rangle &  =\int\left(  \nabla v\cdot\nabla
\phi+V(x)f(v)\dot{f}(v)\phi\right)  \text{,}\\
\langle Q^{\prime\prime}(v)\phi,\psi\rangle &  =\int\left\{  \nabla\phi
\cdot\nabla\psi+V(x)\left[  \dot{f}^{2}(v)+f(v)\ddot{f}(v)\right]  \phi
\psi\right\}
\end{align*}
for all $v,\phi,\psi\in X$. In particular, since $f(0)=0$ and $\dot{f}(0)=1$,
we have $Q^{\prime}(0)=0$ and%
\[
\langle Q^{\prime\prime}(0)\phi,\psi\rangle=\int\left(  \nabla\phi\cdot
\nabla\psi+V(x)\phi\psi\right)  \text{.}%
\]
Now applying the Taylor formula we get%
\begin{align}
Q(v) &  =Q(0)+\langle Q^{\prime}(0),v\rangle+\frac{1}{2!}\langle Q^{\prime
\prime}(0)v,v\rangle+o(\Vert v\Vert ^{2})\nonumber\\
&  =\frac{1}{2}\int\left( \vert \nabla v\vert ^{2}+V(x)v^{2}%
\right)  +o(\Vert v\Vert ^{2})\text{,\qquad as }\Vert
v\Vert \rightarrow0\text{.}\label{f1}%
\end{align}
On the other hand, by the properties of $f$ again, we deduce from \eqref{e2}
that%
\[
\lim_{t\rightarrow0}\frac{G(f(t))}{t^{2}}=\lim_{t\rightarrow0}\frac
{G(f(t))}{f^{2}(t)}\frac{f^{2}(t)}{t^{2}}=0\text{.}%
\]
From this and condition $\left(  g_{0}\right)  $, for any $\varepsilon>0$,
there exists $C_{\varepsilon}>0$ such that%
\[
\vert G(f(t))\vert \leq\varepsilon t^{2}+C_{\varepsilon}\vert
t\vert ^{p/2}\text{.}%
\]
Therefore, since $p>4$, as $\Vert v\Vert \rightarrow0$ we have%
\[
\int G(f(v))=o\left( \Vert v\Vert ^{2}\right)  \text{.}%
\]
Using this and \eqref{f1}, as $\Vert v\Vert \rightarrow0$ we have%
\begin{align*}
\Phi(v) &  =Q(v)-\int G(f(v))\\
&  =\frac{1}{2}\int\left(\vert \nabla v\vert ^{2}+V(x)v^{2}%
\right)  +o(\Vert v\Vert ^{2})\text{.}%
\end{align*}
From this and \eqref{qe}, it is easy to see that the conclusion of the lemma is true.
\end{pf}

\begin{rem}
\label{rk0}If the quadratic form%
\[
\mathfrak{B}(v)=\int\left(  \vert \nabla v\vert ^{2}+V(x)v^{2}%
\right)
\]
is positive definite, then $X^{-}=\left\{  0\right\}  $ and the zero function
$0$ is a strict local minimizer of $\Phi$. Under some additional assumptions
on the nonlinearity $g(u)$, we could deduce that $\Phi$ satisfies the mountain
pass geometry. This approach seems simpler than those presented in
\cite{MR2580150,MR3160456,MR2927594}.
\end{rem}

To apply Proposition \ref{ap1}, in what follows, we shall investigate the
critical groups of $\Phi$ at infinity. Firstly, we prove a crucial property of
the transformation $f$.

\begin{lem}
\label{l0}Suppose that $\tilde{g}:\mathbb{R}\rightarrow\mathbb{R}$ satisfies
$\tilde{g}(s)s\geq0$, then for all $s\in\mathbb{R}$ we have%
\begin{equation}
\tilde{g}(f(s))\dot{f}(s)s\geq\frac{1}{2}\tilde{g}(f(s))f(s)\text{.}
\label{ee}%
\end{equation}

\end{lem}

\begin{pf}
If $s\geq0$ then $\tilde{g}(f(s))\geq0$, by Proposition \ref{p1} \eqref{ite1}
we get \eqref{ee}. If $s<0$, then $\tilde{g}(f(s))\leq0$ and since $\dot{f}$ is
an even function, by Proposition \ref{p1} \eqref{ite1} we have%
\[
-\dot{f}(s)s=\dot{f}(-s)(-s)\geq\frac{f(-s)}{2}\text{.}%
\]
Multiplying both side by $-\tilde{g}(f(s))$ and noting that $f$ is odd, we get%
\[
\tilde{g}(f(s))\dot{f}(s)s\geq-\tilde{g}(f(s))\frac{f(-s)}{2}=\frac{1}%
{2}\tilde{g}(f(s))f(s)\text{.}%
\]

\end{pf}

\begin{lem}
\label{l4}There exists $A>0$ such that, if $\Phi(v)\leq-A$, then%
\[
\left.  \frac{d}{dt}\right\vert _{t=1}\Phi(tv)<0\text{.}%
\]

\end{lem}

\begin{pf}
Otherwise, there exists a sequence $\left\{  v_{n}\right\}  \subset X$ such
that $\Phi(v_{n})\leq-n$ but%
\begin{equation}
\langle\Phi^{\prime}(v_{n}),v_{n}\rangle=\left.  \frac{d}{dt}\right\vert
_{t=1}\Phi(tv_{n})\geq0\text{.}\label{e5}%
\end{equation}
Note that, since $\tilde{g}(s)s\geq0$, by Proposition \ref{p1} \eqref{ite2}
and Lemma \ref{l0}, we have%
\begin{align}
\langle\Phi^{\prime}(v_{n}),v_{n}\rangle &  =\int\left( \vert \nabla
v_{n}\vert ^{2}+\tilde{V}(x)f(v_{n})\dot{f}(v_{n})v_{n}\right)
-\int\tilde{g}(f(v_{n}))\dot{f}(v_{n})v_{n}\nonumber\\
&  \leq\int\left( \vert \nabla v_{n}\vert ^{2}+\tilde{V}%
(x)f^{2}(v_{n})\right)  -\frac{1}{2}\int\tilde{g}(f(v_{n}))f(v_{n}%
)\text{.}\label{e6}%
\end{align}
Combining \eqref{e5}, \eqref{e6} and using \eqref{ge}, we have%
\begin{align}
0 &  \geq\frac{\mu}{2}\Phi(v_{n})-\langle\Phi^{\prime}(v_{n}),v_{n}%
\rangle\nonumber\\
&  \geq\left(  \frac{\mu}{4}-1\right)  \int\left( \vert \nabla
v_{n}\vert ^{2}+\tilde{V}(x)f^{2}(v_{n})\right)  +\frac{1}{2}\int\left(
\tilde{g}(f(v_{n}))f(v_{n})-\mu\tilde{G}(f(v_{n})\right)  \nonumber\\
&  \geq\left(  \frac{\mu}{4}-1\right)  \int\left(\vert \nabla
v_{n}\vert ^{2}+\tilde{V}(x)f^{2}(v_{n})\right)  -\frac{1}{2}\left(
\frac{\mu}{2}-1\right)  m\int f^{2}(v_{n})\text{.}\label{e9}%
\end{align}
We denote%
\[
\rho_{n}=\left\{  \int\left( \vert \nabla v_{n}\vert ^{2}%
+\tilde{V}(x)f^{2}(v_{n})\right)  \right\}  ^{1/2}\text{.}%
\]
Then we must have $\rho_{n}\rightarrow+\infty$. Otherwise, by the argument of
Wu \cite[Page 2626--2628]{MR3160456} mentioned in the proof of Lemma \ref{l1},
$\left\{  v_{n}\right\}  $ is a bounded sequence in $X$ and $\left\{
\Phi(v_{n})\right\}  $ is also bounded, contradicting the fact that
$\Phi(v_{n})\leq-n$.

Set $h_{n}=f(v_{n})/\rho_{n}$. As in \eqref{e7}, $\left\{  h_{n}\right\}  $ is a
bounded sequence in $X$. Up to a subsequence, we may assume that%
\begin{equation}
h_{n}\rightharpoonup h\text{ in }X\text{,\qquad}h_{n}\rightarrow h\text{ in
}L^{2}\text{,\qquad}h_{n}\rightarrow h\text{ a.e. in }\mathbb{R}^{N}\text{.}
\label{e11}%
\end{equation}
Multiplying both sides of \eqref{e9} by $\rho_{n}^{-2}$, we obtain%
\begin{equation}
\frac{1}{2}\left(  \frac{\mu}{2}-1\right)  m\vert h_{n}\vert
_{2}^{2}\geq\frac{\mu}{4}-1\text{.} \label{e10}%
\end{equation}
Then, since $\mu>4$, \eqref{e11} and \eqref{e10} imply that $h\neq0$.

Now, because $\rho_{n}\rightarrow+\infty$, similar to \eqref{e8} we have%
\[
\frac{1}{\rho_{n}^{2}}\int\tilde{g}(f(v_{n}))f(v_{n})\rightarrow
+\infty\text{.}%
\]
Consequently, by \eqref{e5} and \eqref{e6},%
\begin{align*}
0  &  \leq\langle\Phi^{\prime}(v_{n}),v_{n}\rangle\leq\int\left(\vert
\nabla v_{n}\vert ^{2}+\tilde{V}(x)f^{2}(v_{n})\right)  -\frac{1}{2}%
\int\tilde{g}(f(v_{n}))f(v_{n})\\
&  =\rho_{n}^{2}\left(  1-\frac{1}{2\rho_{n}^{2}}\int\tilde{g}(f(v_{n}%
))f(v_{n})\right)  \rightarrow-\infty\text{.}%
\end{align*}
This is impossible. Thus the conclusion of the lemma must be true.
\end{pf}

\begin{lem}
\label{l3}Under our assumptions, $C_{i}(\Phi,\infty)\cong0$ for $i\in
\mathbb{N}$.
\end{lem}

\begin{pf}
Let $B$ be the unit ball in $X$, $S=\partial B$ the unit sphere. Let $A>0$ be
the number given in Lemma \ref{l4}. Without lose of generality we assume that%
\[
-A<\inf_{\Vert v\Vert \leq2}\Phi(v)\text{.}%
\]
Then for $w\in S$, since $\vert f(t)\vert \leq\vert
t\vert $, reasoning as in \eqref{e8} we deduce%
\begin{align*}
\Phi(sw)  &  \leq\frac{1}{2}\int\left(  \vert \nabla(sw)\vert
^{2}+\tilde{V}(x)(sw)^{2}\right)  -\int\tilde{G}(f(sw))\\
&  =s^{2}\left(  \frac{1}{2}-\frac{1}{s^{2}}\int\tilde{G}(f(sw))\right)
\rightarrow-\infty
\end{align*}
as $s\rightarrow+\infty$.

Consequently, there exists $s_{w}>0$ such that $\Phi(s_{w}w)=-A$. Set
$v=s_{w}w$, a direct computation and Lemma \ref{l4} gives%
\[
\left.  \frac{d}{ds}\right\vert _{s=s_{w}}\Phi(sw)=\frac{1}{s_{w}}\left.
\frac{d}{dt}\right\vert _{t=1}\Phi(tv)<0\text{.}%
\]

By the implicit function theorem, $w\mapsto s_{w}$ is a continuous function on
$S$. Using this function and a standard argument (see, e.g.,
\cite{MR1863785,MR1094651}), we can construct a deformation from $X\backslash
B$ to $\Phi_{-A}=\Phi^{-1}(-\infty,-A]$, and deduce via the homotopic
invariance of singular homology%
\[
C_{i}(\Phi,\infty)=H_{i}(X,\Phi_{-A})\cong H_{i}(X,X\backslash
B)=0\text{,\qquad for }i\in\mathbb{N}\text{.}%
\]

\end{pf}

\begin{pf}
[Proof of Theorem \ref{t1}]We have verified that $\Phi$ satisfies the Cerami
condition. By Lemma \ref{l2}, $\Phi$ has a local linking at $0$ with respect to the
decomposition $X=X^{-}\oplus X^{+}$, hence by Proposition \ref{ap2}, for
$\ell=\dim X^{-}$, we have%
\[
C_{\ell}(\Phi,0)\neq0\text{.}%
\]
On the other hand, Lemma \ref{l3} says that $C_{\ell}(\Phi,\infty)=0$. By
Proposition \ref{ap1}, $\Phi$ has a nonzero critical point $v$. Now $u=f(v)$
is a nontrivial solution of Problem \eqref{e1}.
\end{pf}

\section{Multiplicity result}

To prove Theorem \ref{t2}, we will apply the following symmetric mountain pass
theorem due to Ambrosetti-Rabinowitz \cite{MR0370183}.

\begin{prop}
[{\cite[Theorem 9.12]{MR765240}}]\label{smt}Let $X$ be an infinite dimensional
Banach space, $\Phi\in C^{1}(X,\mathbb{R})$ be even, satisfies Cerami
condition and $\Phi(0)=0$. If $X=Y\oplus Z$ with $\dim Y<\infty$, and $\Phi$ satisfies

\begin{enumerate}
\item[$\left(  I_{1}\right)  $] there are constants $\rho,\alpha>0$ such that
$\Phi|_{\partial B_{\rho}\cap Z}\geq\alpha$,

\item[$\left(  I_{2}\right)  $] for any finite dimensional subspace $W\subset
X$, there is an $R=R(W)$ such that $\Phi\leq0$ on $W\backslash B_{R(W)}$,
\end{enumerate}
then $\Phi$ has a sequence of critical values $c_{j}\rightarrow+\infty$.
\end{prop}

\begin{lem}
\label{l5}Let $W$ be a finite dimensional subspace of $X$, then $\Phi$ is
anti-coercive on $W$, that is%
\[
\Phi(v)\rightarrow-\infty\text{,\qquad as }\Vert v\Vert
\rightarrow\infty\text{, }v\in W\text{.}%
\]

\end{lem}

\begin{pf}
For any $\left\{  v_{n}\right\}  \subset W$ with $\Vert v_{n}\Vert
\rightarrow\infty$, set $h_{n}=\Vert v_{n}\Vert ^{-1}v_{n}$. Then
$\left\{  h_{n}\right\}  $ is a bounded sequence in $W$. Because $\dim
W<\infty$, there exists $h\in W\backslash\left\{  0\right\}  $ such that%
\[
\Vert h_{n}-h\Vert \rightarrow0\text{,\qquad\qquad}h_{n}%
(x)\rightarrow h(x)\text{\quad a.e. }x\in\mathbb{R}^{N}\text{.}%
\]
For $x\in\left\{  h\neq0\right\}  $, we have $\vert v_{n}(x)\vert
\rightarrow\infty$ and $\vert f(v_{n}(x))\vert \rightarrow\infty$.
Therefore, if $n$ large enough we have $\vert v_{n}(x)\vert \geq1$.
By Proposition \ref{p1} \eqref{ite3}, we have%
\begin{align*}
\frac{\tilde{G}(f(v_{n}(x)))}{\Vert v_{n}\Vert ^{2}}  &
=\frac{\tilde{G}(f(v_{n}(x)))}{f^{4}(v_{n})}\frac{f^{4}(v_{n}(x))}{v_{n}%
^{2}(x)}h_{n}^{2}(x)\\
& \geq\kappa^{4}\frac{\tilde{G}(f(v_{n}(x)))}{f^{4}(v_{n}(x))}h_{n}%
^{2}(x)\rightarrow+\infty\text{.}%
\end{align*}
By the Fatou lemma, we have%
\[
\int\frac{\tilde{G}(f(v_{n}))}{\Vert v_{n}\Vert ^{2}}\geq
\int_{h\neq0}\frac{\tilde{G}(f(v_{n}))}{\Vert v_{n}\Vert ^{2}%
}\rightarrow+\infty
\]
and%
\begin{align*}
\Phi(v_{n}) &  =\Vert v_{n}\Vert ^{2}\left(  \frac{1}{2\Vert
v_{n}\Vert ^{2}}\int\left(  \vert \nabla v_{n}\vert
^{2}+\tilde{V}(x)f^{2}(v_{n})\right)  -\int\frac{\tilde{G}(f(v_{n}%
))}{\Vert v_{n}\Vert ^{2}}\right)  \\
&  \leq\Vert v_{n}\Vert ^{2}\left(  \frac{1}{2}-\int\frac{\tilde
{G}(f(v_{n}))}{\Vert v_{n}\Vert ^{2}}\right)  \rightarrow
-\infty\text{,}%
\end{align*}
the desired result follows.
\end{pf}

\begin{pf}
[Proof of Theorem \ref{t2}]Under the assumptions of Theorem \ref{t2}, we know
that $\Phi$ is an even functional satisfying the Cerami condition. Lemma
\ref{l5} implies that $\Phi$ satisfies condition $\left(  I_{2}\right)  $ of
Proposition \ref{smt}.

By condition $\left(  g_{0}\right)  $ and Proposition \ref{p1}, there exist
positive constants $C_{1}$ and $C_{2}$ such that%
\begin{equation}
\vert G(f(t))\vert \leq C_{1}\vert t\vert ^{2}%
+C_{2}\vert t\vert ^{p/2}\text{.}\label{G}%
\end{equation}
For $i\geq\ell$, let $Z_{i}=\overline{\operatorname*{span}}\left\{  \phi
_{i},\phi_{i+1},\ldots\right\}  $. Then%
\[
\beta_{i}=\sup_{v\in Z_{i},\Vert v\Vert =1}\vert v\vert
_{2}\rightarrow0\text{,\qquad as }i\rightarrow\infty\text{,}%
\]
see e.g.\ \cite[Lemma 2.5]{MR2548724}. Therefore, there exists $k>\ell$ such that%
\[
\lambda=\eta-C_{1}\beta_{k}^{2}>0\text{.}%
\]
Let%
\[
Y=\operatorname*{span}\left\{  \phi_{1},\ldots,\phi_{k-1}\right\}
\text{,\qquad}Z=\overline{\operatorname*{span}}\left\{  \phi_{k},\phi
_{k+1},\ldots\right\}  \text{.}%
\]
Then $Z\subset X^{+}$ and $X=Y\oplus Z$.

For $v\in Z$, using \eqref{G} and Taylor expansion as in the proof of Lemma
\ref{l2}, and noting that $p>4$, we have%
\begin{align*}
\Phi(v) &  =\frac{1}{2}\int\left(  \vert \nabla v\vert
^{2}+V(x)f^{2}(v)\right)  -\int G(f(v))\\
&  =\frac{1}{2}\int\left(  \vert \nabla v\vert ^{2}+V(x)v^{2}%
\right)  -\int G(f(v))+o(\Vert v\Vert ^{2})\\
&  \geq\eta\Vert v\Vert ^{2}-C_{1}\vert v\vert _{2}%
^{2}-C_{2}\vert v\vert _{p/2}^{p/2}+o(\Vert v\Vert
^{2})\\
&  \geq\eta\Vert v\Vert ^{2}-C_{1}\beta_{k}^{2}\Vert
v\Vert ^{2}+o(\Vert v\Vert ^{2})\\
&  =\lambda\Vert v\Vert ^{2}+o(\Vert v\Vert
^{2})\text{.}%
\end{align*}
Now it is easy to see that condition $\left(  I_{1}\right)  $ of Proposition
\ref{smt} is verified.

By Proposition \ref{smt}, $\Phi$ has a sequence of critical
points $\{v_{n}\}$ such that $\Phi(v_{n})\to+\infty$. Since $\Phi(v)=J(f(v))$,
let $u_{n}=f(v_{n})$. Then $\{u_{n}\}$ is a sequence of solutions for \eqref{e1}
such that $J(u_{n})\to+\infty$.
\end{pf}

\renewcommand{\href}[2]{#2}

\end{document}